\documentstyle[12pt,leqno,amssymb,graphics]{amsart}

\newtheorem{thm}{Theorem}[section]
\newtheorem{lemma}[thm]{Lemma}

\theoremstyle{definition}

\theoremstyle{remark}

\begin{document}

\newcommand{\ct}{\cite}
\newcommand{\pr}{\protect\ref}
\newcommand{\su}{\subseteq}
\newcommand{\pa}{{\partial}}
\newcommand{\im}{{Imm(F,\E)}}
\newcommand{\hf}{{1 \over 2}}
\newcommand{\R}{{\Bbb R}}
\newcommand{\Z}{{\Bbb Z}}
\newcommand{\E}{{{\Bbb R}^3}}
\newcommand{\C}{{{\Bbb Z}/2}}
\newcommand{\lm}{{\lambda}}

\newcommand{\I}{{\mathrm{Id}}}
\newcommand{\4}{{\mathcal{H}}}
\newcommand{\hc}{{H_1(F,\C)}}
\newcommand{\ak}{{ \{ a_k \}  }}   
\newcommand{\bk}{{ \{ b_k \}  }}   
\newcommand{\akk}{{ \{ a'_k \} }}
\newcommand{\bkk}{{ \{ b'_k \} }}

\newcommand{\tb}{{ \Leftrightarrow }} 
\newcommand{\bn}{{ \leftrightarrow }} 

\newcounter{numb}

\title{Complementary regions for immersions of surfaces}
\author{Tahl Nowik}
\address{Department of Mathematics, Bar-Ilan University, 
Ramat-Gan 52900, Israel}
\email{tahl@@math.biu.ac.il}
\date{January 7, 2007}

\begin{abstract}
Let $F$ be a closed surface and $i:F \to S^3$ a generic immersions. 
Then $S^3 - i(F)$ is a union of connected regions, which may be separated into two sets
$\{ U_j \}$ and $\{ V_j \}$ by a checkerboard coloring.
For $k \geq 0$ let $a_k$, $b_k$ be the number of
components $U_j$, $V_j$ with $\chi(U_j) = 1-k$, $\chi(V_j)=1-k$, respectively.
Two more integers attached to $i$ are
the number $N$ of triple points of $i$, and $\chi=\chi(F)$.
In this work we determine what sets of data $(\ak, \bk, \chi, N)$ may appear in this way.
\end{abstract}

\maketitle

\section{The setting, and statement of result}
For $F$ a closed surface and $i: F\to S^3$ a generic immersion, 
we will be interested in the connected components of $S^3 - i(F)$. 
They will be called
the \emph{complementary regions} of $i$, or simply the \emph{regions} of $i$.
Choose one point $p_0 \in S^3 - i(F)$ and color it black. This determines a color black or white for any
point in $S^3 - i(F)$ according to the following prescription: 
If $p \in S^3 - i(F)$ we connect $p$ to $p_0$ with a curve $\gamma$
in general position with respect to $i(F)$, and we color $p$ black or white
according to whether $\gamma$ intersects $i(F)$ an even or odd number of times, respectively. 
This is indeed well defined by simple connectivity of $S^3$.
The color of points within one region is constant and so we may
refer to the color of a region, and this color changes whenever
crossing a sheet of $i(F)$ from one region to a neighboring one. 
We will be interested in the collection of Euler characteristics that may appear for
the set of regions, so we first prove:

\begin{lemma}\label{l1}
If $U$ is a complementary region of $i:F \to S^3$ then $\chi(U) \leq 1$.
\end{lemma}

\begin{pf}
Let $\pa U$ be the natural notion of a boundary for $U$.
It is enough to show that $\pa U$ is connected, since then $\chi(\pa U) \leq 2$
and so $\chi(U) = \hf \chi(\pa U) \leq 1$.
So assume $\pa U$ has at least two connected components $S_1,S_2$, and 
let $T \su U$ be a surface parallel 
to $S_1$. There is a path in $U$ from $S_1$ to $S_2$ crossing $T$ precisely once, 
and since $i(F)$ is connected and disjoint from $T$, this path 
can be completed to a loop in $S^3$ crossing $T$ precisely once, 
contradicting the simple connectivity of $S^3$.
\end{pf}

Given a generic immersion $i:F \to S^3$, color $S^3 - i(F)$ as above, and
we define two sequences $a_0,a_1,a_2,\dots$ and $b_0,b_1,b_2,\dots$ of non-negative integers as follows:
Let $a_k = a_k(i)$ be the number of black regions $U$ with $\chi(U)=1-k$ 
and let $b_k=b_k(i)$ be the number of white regions $U$ with $\chi(U)=1-k$.
We attach two more integers to such immersion, the number $N=N(i)$ of triple points of $i$,
and $\chi = \chi(F)$.
Our goal in this work is to determine what sets of data 
$(\ak,\bk,\chi,N)$ may arise in this way.
We will prove:

\begin{thm}\label{main}
Let $\{ a_k \}_{k \geq 0}$, $\{ b_k \}_{k \geq 0}$ 
be two sequences of non-negative integers which are not identically 0. 
Let $\chi$ and $N$ be integers satisfying $\chi \leq 2$, $N\geq 0$. 
Then there is a closed surface $F$ with $\chi(F) = \chi$ and an immersion
$i:F \to S^3$ with $N$ triple points which realizes the sequences $\ak,\bk$,  iff the 
following equation holds: 
$$\sum_k (1-k)a_k = \sum_k (1-k)b_k = \hf (\chi + N).$$
\end{thm}

The ``only if'' part of Theorem \pr{main} will be proved in Section \pr{onif} and the ``if''
part in Section \pr{if}.

\section{Any immersion satisfies the equation}\label{onif}

In this section we show that any immersion $i:F \to S^3$ satisfies the 
equalities of Theorem \pr{main}. This has been shown for orientable surfaces
in \ct{f} using an elaborate order 1 invariant appearing there. 
We will present a self contained proof which
covers both orientable and non-orientable surfaces. 
For generic immersion $i: F \to S^3$ let $A(i)$ be the union of all black regions and $B(i)$ 
the union of all white regions of $i$, then 
$\sum_k (1-k)a_k = \chi(A(i))$ and $\sum_k (1-k)b_k = \chi(B(i))$.
We now define $a(i)$ and $b(i)$ by $a(i) = \chi(A(i)) - \hf N(i)$ and $b(i) = \chi(B(i)) - \hf N(i)$ 
(this may be a half integer). Our goal is then to show that for any immersion $i:F \to S^3$, 
$a(i) = b(i) = \hf \chi(F)$. We will show this by first showing 
that $a(i)$ and $b(i)$ are 
constant along regular homotopies, and then verifying the equalities for one chosen 
immersion in every regular homotopy class. 

Along a generic regular homotopy, there are some finitely many 
times when the topology of the complementary regions may change, namely, at the times when the regular homotopy 
passes through non-stable immersion. There are four types of such occurrences, which we name $E,H,T,Q$. 
They appear in Figure 1, and a model immersion 
in $\E$ for each configuration is given below, depending on a parameter $\lm$. The immersions depicted in Figure 1 
correspond to some $\lm > 0$, and the unstable immersion corresponds to $\lm = 0$. 
We will now see that for each of the four types of occurrences, 
if we continuously carry the coloring along the regular homotopy,
by keeping it fixed on points which are remote from the given occurrence, 
then the topology of the regions and the number of triple points 
change in a way such that $a(i),b(i)$ remain constant.

\begin{figure}[t]
\scalebox{0.7}{\includegraphics{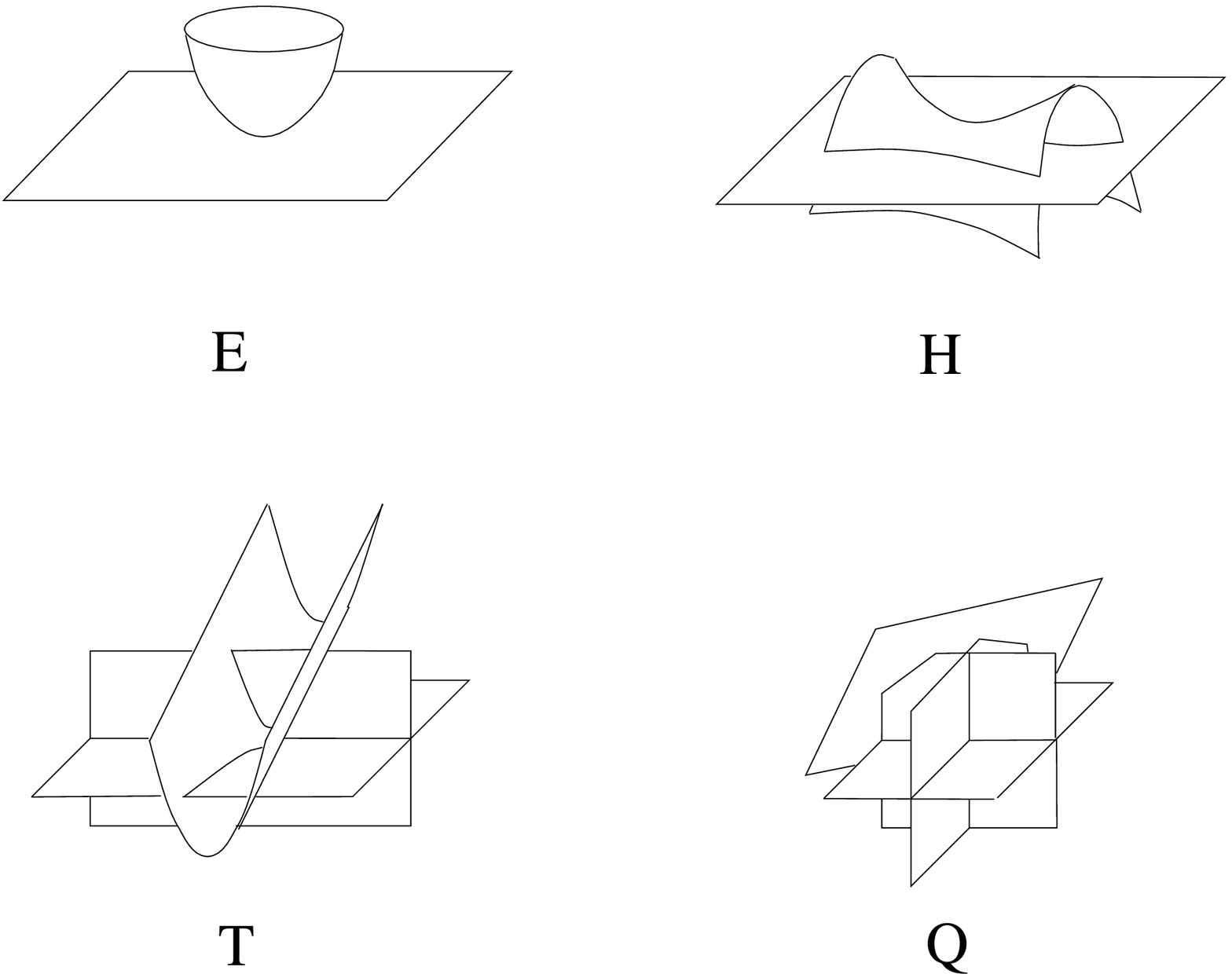}}
\caption{}\label{f1}
\end{figure}

Model for $E$: \ $z=0$, \ $z=x^2+y^2+\lm$.
When moving from $\lm > 0$ to $\lm < 0$ a new 3-cell region appears. Say this region is black
then it adds 1 to $a(i)$. The second change occurring, is that the region just above
the plane $z=0$ has a 2-handle removed from it.
This region is seen to be black, since one needs to cross 
both sheets involved in the configuration in order to pass from 
this region into the new 3-cell. 
And so all together $a(i)$ remains fixed.
The topology of the white regions has remained unchanged and so 
$b(i)$ remains fixed as well.

Model for $H$: \  $z=0$,  \ $z=x^2-y^2+\lm$. 
The changes occurring in the 
neighboring regions when passing from $\lm > 0$ to $\lm < 0$ are that a 1-handle is removed
from the region $X$ just above the $x$ axis, and a 1-handle is added to the region $Y$
just below the $y$ axis. It is seen that the color of these two regions is the same 
since they are on opposite sides of both sheets involved in the configuration, and so $a(i),b(i)$ remain fixed.

Model for $T$:  \ $z=0$, \ $y=0$,  \ $z=y+x^2+\lm$.
When passing from $\lm > 0$ to $\lm < 0$ a new 3-cell region appears.
The second change occurring is that a 1-handle is removed from the region near the
$x$ axis having negative $y$ values and positive $z$ values. The color of this region is seen to be opposite
to that of the new 3-cell, since it is on the other side of all three sheets involved in the configuration, 
and so $\chi(A(i))$ and $\chi(B(i))$ both increase by 1. 
But on the other hand, 
two triple points are added to $i$ during this occurrence, 
and so all together both $a(i)$ and $b(i)$ remained fixed.

Model for $Q$:  \ $z=0$, \ $y=0$, \ $x=0$, \ $z=x+y+\lm$. 
When passing from $\lm >0$ to $\lm < 0$, a simplex region vanishes, 
and a new simplex region appears. They are of the same color since they are
on the opposite side of all four sheets involved in the configuration.
All other regions are topologically unchanged and so again, $a(i)$ and $b(i)$ 
remain fixed.

We have thus shown that $a(i)$ and $b(i)$ are constant on every regular homotopy class.
It remains to find this 
value by direct inspection of one immersion in every regular homotopy class,
and see that indeed $a(i)=b(i)=\hf \chi(F)$. It is shown in \ct{p} that if $F$ is orientable then any
immersion $i:F \to \E$ is regularly homotopic to either a standard embedding, or a standard embedding with
an additional ``ring'' added along a particular circle in the surface.
(A ring along a trivial circle is seen in Figure 3b). 
It follows that the same is true for immersions in $S^3$.
For a standard embedding the regions on the two sides
indeed have $\chi(U) = \hf \chi(F)$. If a ring is added, this adds a solid torus region, whose $\chi$ is 0,
and it does not change the topology of the previously existing regions,
and so the equalities still hold. It is further shown in \ct{p} that if $F$ is non-orientable, then any immersion
$i:F \to \E$ is regularly homotopic to a connect sum of Boy's surfaces (right and left handed), 
where by connect sum of immersions we mean the following. If $i:F \to S^3$, is an immersion
of a closed surface into $S^3$, and if $p \in i(F)$ is a point not in the intersection set of $i$, 
then let $B$ be a small 3-ball neighborhood of $p$ disjoint from the intersection set of $i$ and 
delete $B$ from $S^3$, obtaining an immersion $i|_{F-D}:F-D \to S^3 - B$  where $D = i^{-1}(B)$
is a disc in $F$.
If $i':F' \to S^3$ is another such immersion, do the same with some $B',D'$. Now glue $S^3-B$ to $S^3-B'$
along their boundaries,
so that the boundaries of $F-D$ and $F'-D'$ will match, obtaining an immersion $F \# F' \to S^3 \# S^3 =S^3$ which we 
call a connect sum of $i,i'$.
When such connect sum operation is performed, then the two complementary regions of $i$ on the two sides of $i(D)$,
merge with the corresponding complementary regions of $i'$ on the two sides of $i'(D')$.
Each such merger is along a disc, namely, a hemisphere of $\pa B$.

In \ct{p} Boy's surface is depicted with a small window removed, so one can peak inside and convince 
one's self that the complementary regions of Boy's surface are two 3-cells, one of each color. 
Also, Boy's surface has one triple point. And so
the complementary regions of a connect sum of $N$ Boy's surfaces 
are still two 3-cells, one of each color, and 
the number of triple points is $N$. And so $a(i)=b(i)=1-\hf N = \hf\chi(F)$.
This completes the proof of the ``only if'' part of Theorem \pr{main}.

We make the following remark. The ``only if'' part of Theorem \pr{main} 
is a refinement 
of what may be deduced in this setting from the 
result of Izumiya and Marar (\ct{im}), namely $\chi(i(F)) = \chi(F) + N$. 
Indeed, present $S^3$ as the union of a regular neighborhood $U$ of $i(F)$, 
and the union $V$ of all complementary regions, slightly diminished, 
so $\chi(V) = \sum_k (1-k)a_k + \sum_k (1-k)b_k$. 
We get $\chi(S^3) = 0 = \chi(i(F)) + \chi(V) - \chi( \pa V) = \chi(i(F)) - \chi(V)$, or $\chi(i(F)) = \chi(V)$. 
This shows that the Izumiya Marar equality is equivalent in this setting to
$\chi(V) = \chi(F)+N$, which is the sum of our two equalities
$\sum_k (1-k)a_k = \hf (\chi + N)$ and
$\sum_k (1-k)b_k = \hf (\chi + N)$.

\section{Realizing the data by immersions}\label{if}

In this section we show that any data $(\ak,\bk,\chi,N)$ 
satisfying the conditions of Theorem \pr{main}, 
may be realized by some $F$ and $i$.
The proof is by induction on $N + \sum_k a_k + \sum_k b_k$.

Let $(\ak,\bk,\chi,N)$ be data satisfying the conditions of Theorem \pr{main}.
If $a_0 \geq 1$ and there is $r \geq 1$ with $a_r \geq 1$, 
then let $\akk$ be the sequence obtained from $\ak$ by 
subtracting 1 from $a_0$, subtracting 1 from $a_r$, and then adding 1 to $a_{r-1}$. 
(so if $r=1$ then $a'_0=a_0$).
If $\ak$ satisfies the conditions of Theorem \pr{main}, 
then so does $\{ a'_k \}$, and we have $a'_{r-1} \geq 1$. 
By induction there is a surface $F$ and immersion $i:F \to S^3$ realizing 
$(\{a'_k\},\bk,\chi,N)$.
Let $U$ be a black complementary region for $i$ with $\chi(U)=1-(r-1)$, and change $i$ in a disc on the boundary of $U$
(which is disjoint from the intersection set of $i$) as in Figure 2. Then
$\chi(U)$ decreases by 1, and there is a new black 3-cell region appearing. The topology of all other regions 
is unchanged, and so the new immersion realizes $(\ak,\bk,\chi,N)$. So we may assume that either (i) $a_0=0$,
or (ii) $a_k=0$ for all $k \geq 1$. By the same argument, the same may be assumed for $\bk$.
If say $\ak$ satisfies (i) and $\bk$ satisfies (ii), then $\sum_k (1-k)a_k \leq  0$, 
and since the sequences are assumed not to be identically 0, $\sum_k (1-k)b_k > 0$, and so they cannot be equal. 
It follows that we may assume either both $\ak,\bk$ satisfy (i) or both satisfy (ii).

\begin{figure}[t]
\scalebox{0.7}{\includegraphics{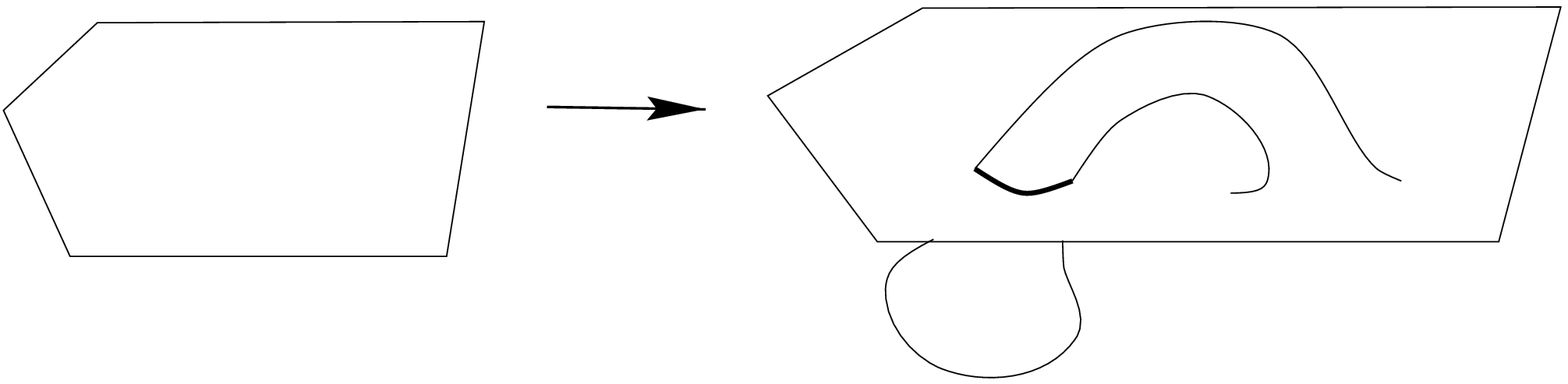}}
\caption{}\label{f2}
\end{figure}

\begin{figure}[t]
\scalebox{0.7}{\includegraphics{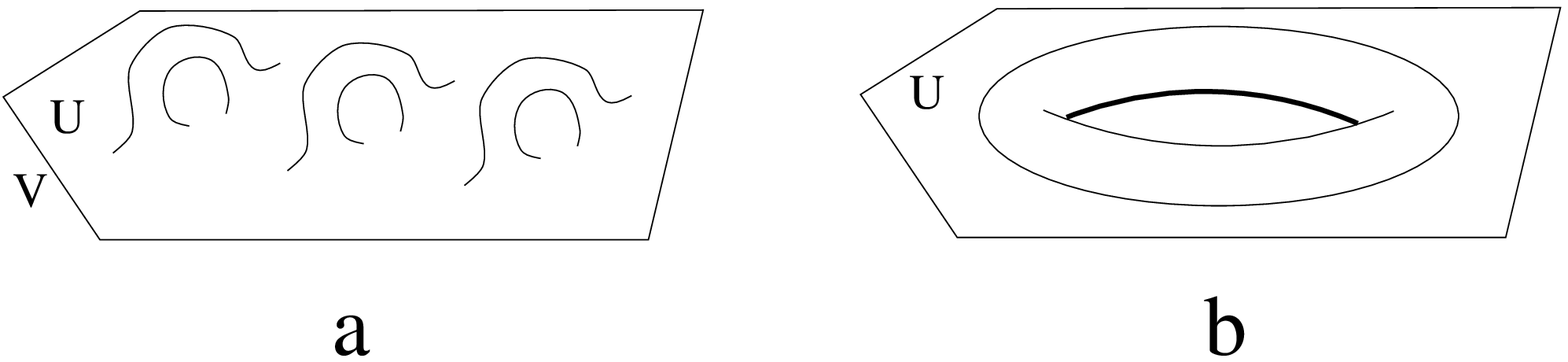}}
\caption{}\label{f3}
\end{figure}

We will first assume both sequences satisfy (i), that is $a_0=b_0=0$, and so 
$\sum_k (1-k)a_k \leq 0$, and so necessarily $\chi \leq 0$. 
If $N \geq 1$,
then by induction we have a surface and immersion $i$ realizing 
$(\ak,\bk,\chi+1,N-1)$.
Take the connect sum of $i$ with Boy's surface, then since Boy's surface has 1 triple point 
and one complementary 3-cell of each color, 
the new immersion realizes $(\ak,\bk,\chi,N)$.
So we may assume $N=0$.

Before dealing with this case, we introduce two more operations on immersions.
Let $i$ be an immersion and let $D \su i(F)$ be a disc disjoint from the intersection set of $i$.
$D$ is part of the boundary of two neighboring 
regions $U,V$ (of opposite color). 
For $g \geq 0$, a $g$-operation on $D$ is an operation as in Figure 3a, 
adding $g$ handles to $F$,
by this reducing $\chi(F)$ by $2g$, and reducing $\chi(U),\chi(V)$ each by $g$.
We will say that this $g$-operation was performed on the pair $U,V$.
The second operation will be called a ring operation, which is adding a \emph{ring} on $D$, say on the side of $U$, as in Figure 3b. 
(Only the immersion changes here, not the topology of the surface.)
The only effect of a ring operation is the creation of
a new solid torus region which is a neighbor of $U$, 
and so of color opposite that of $U$.
The topology of all previously existing regions remains unchanged. 
We will say that this ring operation was performed in the region $U$.

Back to our sequences $\ak,\bk$, recall that we are now assuming $a_0=b_0=0$, $N=0$. 
Assume first that the following holds: There is a $p$ such that $a_k=b_k=0$
for all $k \neq p$, and so (since the sequences are not identically 0) 
$a_p \geq 1$, $b_p \geq 1$. If $p=1$, then by our equalities, $\chi = 0$.
Start with a standard embedding of a torus, and perform 
$a_p - 1$ ring operations in the white side and $b_p - 1$ ring operations in the black side, obtaining
an immersion realizing the data. So we may assume $p \geq 2$ and so necessarily $a_p = b_p$. 
If $a_p=b_p=1$, then $\chi = 2(1-p)$ and we may realize our data by a standard  embedding of an orientable surface 
$F$ of genus $p$. Otherwise $a_p = b_p \geq 2$.
Let $\{a'_k \}$, $\{ b'_k \}$ 
be the sequences  with $a'_p = b'_p = a_p -1$ and $a_k =b_k = 0$ for all $k \neq p$,
and note that necessarily $\chi \leq 2(1-p)$.
By induction we have an immersion $i$ realizing $(\akk,\bkk,\chi+2(p-1),0)$. 
Now perform a ring operation on $i$
in some arbitrary place, creating a new solid torus region $U$. In $U$ perform another ring operation
creating another solid torus region $V$ which is a neighbor of $U$ (and so of
opposite color).
Now perform a $(p-1)$-operation on the pair $U,V$ to obtain an 
immersion realizing $(\ak,\bk,\chi,0)$. 

If there is no $p$ as above, then there are some 
$r\neq s$ such that $a_r \geq 1$ and $b_s \geq 1$, and assume $r < s$. Let
$\akk,\bkk$ be the sequences obtained from $\ak,\bk$ in the following way. Subtract 1 from $a_r$,
subtract 1 from $b_s$ and then add 1 to $b_{s-r+1}$ (so if $r=1$ then the sequence $\bk$ remains unchanged.)
The sequence $\akk$ is not identically 0 since that would imply that 
$\sum_k (1-k)a_k = 1-r$ which could not equal
$\sum_k (1-k)b_k \leq 1-s$. Note also that necessarily $\chi \leq 2(1-r)$.
By induction we may realize $(\akk,\bkk,\chi +2(r-1), 0)$. 
We have $b'_{s-r+1} \geq 1$ and let $U$ be a corresponding region, 
i.e. $U$ is a white region with $\chi(U) = r-s$.
Now, first perform a ring operation in $U$, creating a new black
solid torus region $V$ which is a neighbor
of $U$. Then
perform an $(r-1)$-operation on the pair $U,V$ of neighboring regions, by this decreasing 
$\chi(U)$ and $\chi(V)$ each by $r-1$, and decreasing $\chi(F)$ by $2(r-1)$, and so realizing our original
data $(\ak,\bk,\chi,0)$. (If $r=1$ then the $(r-1)$-operation means doing nothing.)

We are left with the case both sequences satisfy (ii), that is $a_k=b_k=0$ for all $k \geq 1$.
Since the sequences are not identically 0, $\sum_k (1-k)a_k \geq 1$ and so $\chi + N \geq 2$.
If $\chi < 2$ then $N > 0$, and as before we may use a realization of $(\ak,\bk,\chi+1,N-1)$ 
to produce a realization
of $(\ak,\bk,\chi,N)$, so we may assume $\chi =2$ i.e. the surface $F$ is $S^2$. So now $\hf(2+N)$ 
is a sum of integers and so $N$ is even. Our task then, 
is to construct for any even $N$, an immersion
$i:S^2 \to S^3$ with complementary regions which are $\hf (2+N)$ 3-cells of each color. 
From the ``only if'' part of Theorem \pr{main} we know that 
if all complementary regions of an immersion are 3-cells, then necessarily half of 
them are of each color. And so when constructing our immersions below,
we need only verify that the \emph{total} number of 
3-cells is $2+N$.

\begin{figure}[t]
\scalebox{0.7}{\includegraphics{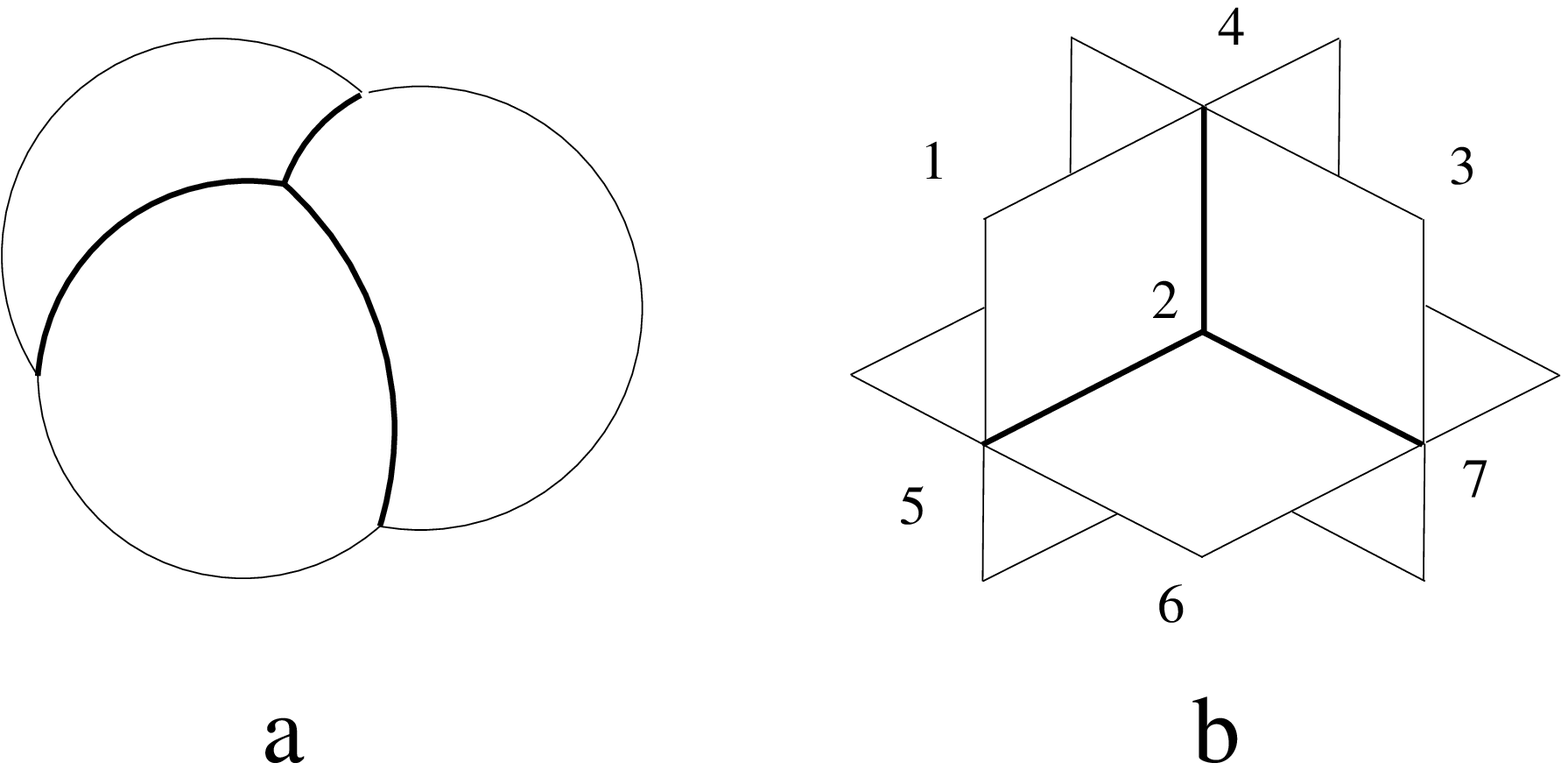}}
\caption{}\label{f4}
\end{figure}

\begin{figure}[t]
\scalebox{0.7}{\includegraphics{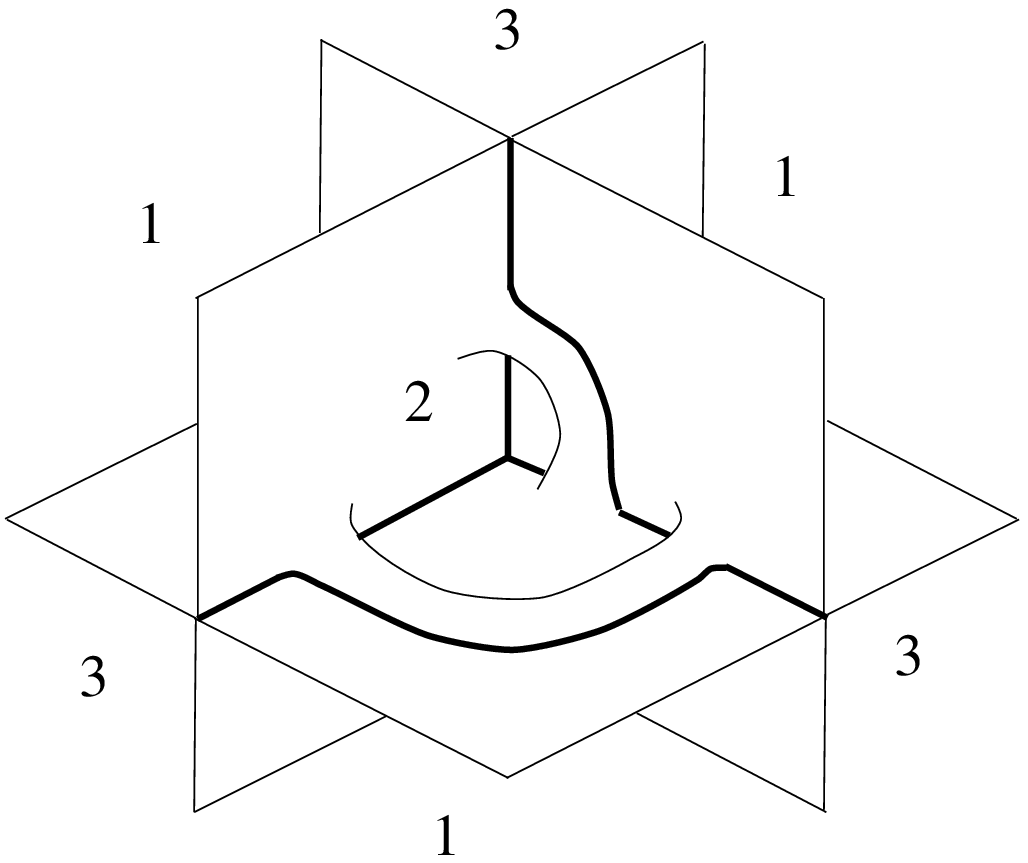}}
\caption{}\label{f5}
\end{figure}

For $N=0$ take an embedding which indeed gives two 3-cell.
For $N=2$, we must construct an immersion $i:S^2 \to S^3$ with 2 triple points and 4 complementary 3-cells.
Start with three spheres, each embedded,
and intersecting each other with two triple points, as in Figure 4a. 
They divide $S^3$ into eight 3-cells.
A small neighborhood of one of the two triple points is 
as appears in Figure 4b.
The 8 different regions appearing
in this neighborhood, belong to the 8 different
complementary regions. In the figure, one of these regions is hidden.

We now attach two tubes to the three spheres, 
merging them into one sphere. We add the tubes as in Figure 5
which we now explain. 
The horizontal tube connects the two vertical sheets, 
and is half way above and half way below the horizontal sheet. 
The upper half of this tube creates a path merging regions 1 and 3 of Figure
4b, and its lower half merges regions 5 and 7 of that figure.
The vertical tube connects the horizontal sheet with the left hand vertical
sheet and is 
half to the left and 
half to the right of the right hand vertical sheet. 
It opens paths merging regions 1 and 6 of Figure 4b, and regions 4 and 7 of that figure.
All together, regions 1,3,6 of Figure 4b have merged into the one 3-cell
region 1 of Figure 5, 
and regions 4,5,7 of Figure 4b have merged into the one 3-cell 
region 3 of Figure 5. And so we have four complementary 3-cells as needed.
This completes the case $N=2$.

For general even $N$, take a connect sum of $\hf N$ copies of the immersion we have constructed for the case $N=2$.
Each additional such immersion adds two new triple points as needed,
and two of its four 3-cell regions merge along discs with two of the existing ones, 
and so precisely two new 3-cell regions are added, as needed.

\end{document}